\documentclass[12pt]{amsart}

\usepackage{amsmath}
\usepackage{amssymb}
\usepackage{url}
\usepackage{amscd}
\usepackage{mathrsfs}
\usepackage[dvips]{graphicx}
\usepackage{indentfirst}
\usepackage{setspace}

\theoremstyle{plain}
\newtheorem{thm}{Theorem}

\newcommand{\gf}{{\mathrm{GF}}}

\newcommand{\cc}{{\mathbf c}}

\newcommand{\tr}{{\mathrm{Tr}}}

\newcommand{\F}{{\mathcal F}}
\newcommand{\A}{{\mathcal A}}

\newcommand{\C}{{\mathcal C}}

\renewcommand{\vec}[1]{\underline{#1}}

\def\({\left(}
\def\){\right)}

\begin{document}

\title{The weight distributions of a class of cyclic codes III}

\author[M. Xiong]{\sc Maosheng Xiong}

\address{Maosheng Xiong: Department of Mathematics,
Hong Kong University of Science and Technology,
Clear Water Bay, Kowloon, Hong Kong
}
\email{mamsxiong@ust.hk}

\thanks{The author was supported by the Research Grants Council of Hong Kong under Project Nos. RGC606211 and DAG11SC02.}


\begin{abstract}
Recently, the weight distributions of the duals of the cyclic codes with two zeros have been obtained for several cases in \cite{DL1,DL2,WT,X,X2}. In this paper we solve one more special case. The problem of finding the weight distribution is transformed into a problem of evaluating certain character sums over finite fields, which in turn can be solved by using the Jacobi sums directly.
\end{abstract}

\maketitle

\thispagestyle{empty}

\maketitle

\thispagestyle{empty}

\section{Introduction}
Let $\gf(q)$ be the finite field of order $q$, where $q=p^s$, $s$ is a positive integer and $p$ is a prime number. An $[n,k,d]$-cyclic code $\C$ over $\gf(q)$ is a $k$-dimensional linear subspace of $\gf(q)^n$ with minimum distance $d$, satisfying the condition that if $(c_0,c_1,\ldots,c_{n-2},c_{n-1}) \in \C$, then the cyclic shift $(c_{n-1},c_0,c_1,\ldots,c_{n-2})$ is also in $\C$. Let $A_i$ denote the number of codewords with Hamming weight $i$ in $\C$. The weight enumerator of $\C$ is defined by
\[1+A_1x+A_2x^2+\cdots+A_nx^n.\]
The sequence $(1,A_1,\ldots,A_n)$ is called the weight distribution of $\C$. In coding theory it is often desirable to know the weight distribution of a code because they contain a lot of important information, for example, they can be used to estimate the error correcting capability and the error probability of error detection and correction with respect to some algorithms. This is quite useful in practice. Many important families of cyclic codes have been studied extensively in the literature, so are their various properties. However the weight distributions are difficult to obtain in general and they are known only for a few special families.

Given a positive integer $m$, let $r=q^m$, and $\alpha$ be a generator of the multiplicative group $\gf(r)^*:=\gf(r)-\{0\}$. Let $h$ be a factor of $q-1$, and $e$ be a factor of $h$. Assume that $h \ge e >1$. Define
\begin{eqnarray} \label{1:para} g=\alpha^{(q-1)/h}, \quad n=\frac{h(r-1)}{q-1},\quad \beta=\alpha^{(r-1)/e},\quad N=\gcd\left(m,\frac{e(q-1)}{h}\right).\end{eqnarray}
It is known that the order of $g$ is $n$, $(g \beta)^n=1$ and the minimal polynomials of $g^{-1}$ and $(\beta g)^{-1}$ are distinct over $\gf(q)$, hence their product is a factor of $x^n-1$ (see \cite{DL1}).
Define the cyclic code over $\gf(q)$ by
\begin{eqnarray} \label{1:code} \C_{(q,m,h,e)}=\left\{\cc_{(a,b)}: a,b \in \gf(r)\right\},\end{eqnarray}
where the codeword $\cc_{(a,b)}$ is give by
\begin{eqnarray} \label{1:codedef} \cc_{(a,b)}:=\left(\tr\left(ag^i+b(\beta g)^i\right)\right)_{i=0}^{n-1}.\end{eqnarray}
Here for simplicity $\tr$ is the trace function from $\gf(r)$ to $\gf(q)$.

The code $\C_{(q,m,h,e)}$ has been an interesting subject of study for a long time. For example, when $h=q-1$, the code $\C_{(q,m,h,e)}$ is the dual of the primitive cyclic linear code with two zeros, which have been studied extensively (see for example \cite{BM, CCD, CCZ,C,Mc,MR,S,YCD}). In general the dimension of $\C_{(q,m,h,e)}$ is $2m$, but determining the weight distribution is very difficult. However, in certain special cases the weight distribution is known. We summarize these cases below.
\begin{itemize}
\item[1)] $e>1$ and $N=1$ (\cite{DL1});

\item[2)] $e=2$ and $N=2$ (\cite{DL1});

\item[3)] $e=2$ and $N=3$ (\cite{DL2});

\item[4)] $e=2$ and $p^j+1 \equiv 0 \pmod{N}$, where $j$ is a positive integer (\cite{DL2});

\item[5)] $e=3$ and $N=2$ (\cite{WT});

\item[6)] $e=4$ and $N=2$ (\cite{X});

\item[7)] $e=3$ and $N=3$ (\cite{X2}).

\end{itemize}

In this paper we compute the weight distribution for $e=3$ under the condition that $p^j +1\equiv 0 \pmod{N}$ for some positive integer $j$. This extends the work \cite{DL2}. The results are as follows.

\begin{thm} \label{thm1} Let $\C_{(q,m,h,e)}$ be the cyclic code defined by (\ref{1:code}) and (\ref{1:codedef}), and the parameters are given by (\ref{1:para}), hence $r=p^{sm}$. Assume that $e=3$, $N \ge 2$ and there exists a positive integer $j$ such that $p^j \equiv -1 \pmod{N}$. Let $j$ be the least such. Assume that $sm=2j\gamma$.
\begin{itemize}
\item[(1).] If $\gamma, p$ and $(p^j+1)/N$ are all odd, then

\begin{itemize}
\item[(1.1)] if $N|\frac{q-1}{h}$, the weight distribution of $\C_{(q,m,h,e)}$ is given by Table \ref{1:t11};

\item[(1.2)] if $N \dagger \frac{q-1}{h}$, the weight distribution of     $\C_{(q,m,h,e)}$ is given by Table \ref{1:t12}.
\end{itemize}

\item[(2).] if one of $\gamma, p$ or $(p^j+1)/N$ is even, then

\begin{itemize}
\item[(2.1)] if $N|\frac{q-1}{h}$, the weight distribution of $\C_{(q,m,h,e)}$ is given by Table \ref{1:t21};

\item[(2.2)] if $N \dagger \frac{q-1}{h}$, the weight distribution of     $\C_{(q,m,h,e)}$ is given by Table \ref{1:t22}.
\end{itemize}

\end{itemize}
\end{thm}

Theorem \ref{thm1} can be compared with results in \cite{WT} and \cite{X2}. In fact, when $e=3,N=2$, the weight distribution of the code was listed in \cite[Table 1]{WT}. On the other hand, for $e=3,N=2$, since $N|(q-1)$, $q$ and $p$ are always odd, so $p \equiv -1 \pmod{N}$, moreover $N|\frac{q-1}{h}$ since $\gcd(e,N)=1$, so the weight distribution is also provided by Table 1 in (1.1) of Theorem \ref{thm1} (or by Table 3 in (2.1), which turns out to be the same.). Actually when $N=2$, Table \ref{1:t11} of Theorem \ref{thm1} matches \cite[Table 1]{WT}. When $e=N=3$, the weight distribution of the code was obtained recently in \cite{X2}. On the other hand, for $e=N=3$, if $p \equiv -1 \pmod{3}$, then at least one of $p,(p+1)/N$ is even, so (2) of Theorem 1 also applies. It turns out that Table \ref{1:t21} of Theorem \ref{thm1} matches \cite[Table 1]{X2} and Table \ref{1:t22} of Theorem \ref{thm1} matches \cite[Table 2]{X2}. For $e=N=3$ and $p \equiv 1 \pmod{3}$ on which Theorem \ref{thm1} does not apply, the weight distribution turns out to be much more complicated (see \cite[Tables 3 and 4]{X2}), with 13 distinct none-zero weights in general. Interested readers may refer to \cite{WT,X2} for some examples of codes that are computed from Magma. It might be interesting to compute examples of codes with $e=3,N \ge 4$ and compare them with Theorem \ref{thm1}. However that seems quite difficult because the parameters of the codes are too large.

It has been known from \cite{X,X2} that to find the weight distribution of the codes, it suffices to evaluate certain character sums over finite fields. In \cite{X,X2} this was succeeded by relating the character sums to counting the number of points on some elliptic curves. In this paper, we use a more direct approach, that is, for $e=3$, the character sums are naturally related to the Jacobi sums, which can be evaluated for the special cases assumed in Theorem \ref{thm1}.


\begin{table}[ht] 
\caption{The case for (1.1) of Theorem \ref{thm1}}
\centering
\begin{tabular}{|c|| c|}
\hline
\hline
Weight & Frequency  \\
[0.5ex]
\hline
$\frac{h}{q}\left\{r-\sqrt{r}(N-1)\right\}$ & $\frac{r-1}{N^3}\bigl\{r+\sqrt{r}(N^2-3N+2)-3N+1\bigr\}$ \\
\hline
$\frac{h}{q}\left\{r+\sqrt{r}\right\}$ & $\frac{(r-1)(N-1)}{N^3}\bigl\{r(N-1)^2-\sqrt{r}(N-2)-(N-1)(2N+1)\bigr\}$ \\
\hline
$\frac{h}{3q}\left\{3r-\sqrt{r}(N-3)\right\}$ & $\frac{3(\sqrt{r}+1)(r-1)(N-1)}{N^3}\bigl\{\sqrt{r}(N-1)-1\bigr\}$ \\
\hline
$\frac{h}{3q}\left\{3r-\sqrt{r}(2N-3)\right\}$ &  $\frac{3(\sqrt{r}+1)(r-1)(N-1)}{N^3}\bigl\{\sqrt{r}-N+1\bigr\}$ \\
\hline
$\frac{2h}{3q}\left\{r-\sqrt{r}(N-1)\right\}$ & $\frac{3(r-1)}{N}$ \\
\hline
$\frac{2h}{3q}\left\{r+\sqrt{r}\right\}$ & $\frac{3(r-1)(N-1)}{N}$ \\
\hline
$0$ & $1$ \\
\hline
\hline
\end{tabular}
\label{1:t11}
\end{table}

\begin{table}[ht] 
\caption{The case for (1.2) of Theorem \ref{thm1}}
\centering
\begin{tabular}{|c|| c|}
\hline
\hline
Weight & Frequency  \\
[0.5ex]
\hline
$\frac{h}{q}\left\{r-\sqrt{r}(N-1)\right\}$ & $\frac{(r-1)(\sqrt{r}+1)^2}{N^3}$ \\
\hline
$\frac{h}{q}\left\{r+\sqrt{r}\right\}$ & $\frac{r-1}{N^3}\bigl\{r(N-1)^3-2\sqrt{r}-(N-1)(2N^2-4N-1)\bigr\}$ \\
\hline
$\frac{h}{3q}\left\{3r-\sqrt{r}(N-3)\right\}$ & $\frac{3(r-1)}{N^3}\bigl\{r(N-1)^2+2\sqrt{r}-2N^2+2N+1\bigr\}$ \\
\hline
$\frac{h}{3q}\left\{3r-\sqrt{r}(2N-3)\right\}$ &  $\frac{3(\sqrt{r}+1)(r-1)}{N^3}\bigl\{\sqrt{r}(N-1)-N-1\bigr\}$ \\
\hline
$\frac{h}{3q}\left\{2r-\sqrt{r}(N-2)\right\}$ & $\frac{6(r-1)}{N}$ \\
\hline
$\frac{2h}{3q}\left\{r+\sqrt{r}\right\}$ & $\frac{3(r-1)(N-2)}{N}$ \\
\hline
$0$ & $1$ \\
\hline
\hline
\end{tabular}
\label{1:t12}
\end{table}

\begin{table}[ht] 
\caption{The case for (2.1) of Theorem \ref{thm1}}
\centering
\begin{tabular}{|c|| c|}
\hline
\hline
Weight & Frequency  \\
[0.5ex]
\hline
$\frac{h}{q}\left\{r+\sqrt{r}(-1)^{\gamma}(N-1)\right\}$ & $\frac{r-1}{N^3}\bigl\{r-\sqrt{r}(-1)^{\gamma}(N^2-3N+2)-3N+1\bigr\}$ \\
\hline
$\frac{h}{q}\left\{r-\sqrt{r}(-1)^{\gamma}\right\}$ & $\frac{(r-1)(N-1)}{N^3}\bigl\{r(N-1)^2+\sqrt{r}(-1)^{\gamma}(N-2)-(N-1)(2N+1)\bigr\}$ \\
\hline
$\frac{h}{3q}\left\{3r+\sqrt{r}(-1)^{\gamma}(N-3)\right\}$ & $\frac{3(r-1)(N-1)}{N^3}\bigl\{r(N-1)-\sqrt{r}(-1)^{\gamma}(N-2)-1\bigr\}$ \\
\hline
$\frac{h}{3q}\left\{3r+\sqrt{r}(-1)^{\gamma}(2N-3)\right\}$ &  $\frac{3(r-1)(N-1)}{N^3}\bigl\{r+\sqrt{r}(-1)^{\gamma}(N-2)-N+1\bigr\}$ \\
\hline
$\frac{2h}{3q}\left\{r+\sqrt{r}(-1)^{\gamma}(N-1)\right\}$ & $\frac{3(r-1)}{N}$ \\
\hline
$\frac{2h}{3q}\left\{r-\sqrt{r}(-1)^{\gamma}\right\}$ & $\frac{3(r-1)(N-1)}{N}$ \\
\hline
$0$ & $1$ \\
\hline
\hline
\end{tabular}
\label{1:t21}
\end{table}

\begin{table}[ht] 
\caption{The case for (2.2) of Theorem \ref{thm1}}
\centering
\begin{tabular}{|c|| c|}
\hline
\hline
Weight & Frequency  \\
[0.5ex]
\hline
$\frac{h}{q}\left\{r+\sqrt{r}(-1)^{\gamma}(N-1)\right\}$ & $\frac{r-1}{N^3}\bigl\{r-2\sqrt{r}(-1)^{\gamma}+1\bigr\}$ \\
\hline
$\frac{h}{q}\left\{r-\sqrt{r}(-1)^{\gamma}\right\}$ & $\frac{r-1}{N^3}\bigl\{r(N-1)^3+2\sqrt{r}(-1)^{\gamma}-(N-1)(2N^2-4N-1)\bigr\}$ \\
\hline
$\frac{h}{3q}\left\{3r+\sqrt{r}(-1)^{\gamma}(N-3)\right\}$ & $\frac{3(r-1)}{N^3}\bigl\{r(N-1)^2-2\sqrt{r}(-1)^{\gamma}-2N^2+2N+1\bigr\}$ \\
\hline
$\frac{h}{3q}\left\{3r+\sqrt{r}(-1)^{\gamma}(2N-3)\right\}$ &  $\frac{3(r-1)}{N^3}\bigl\{r(N-1)+2\sqrt{r}(-1)^{\gamma}-N-1\bigr\}$ \\
\hline
$\frac{h}{3q}\left\{2r+\sqrt{r}(-1)^{\gamma}(N-2)\right\}$ & $\frac{6(r-1)}{N}$ \\
\hline
$\frac{2h}{3q}\left\{r-\sqrt{r}(-1)^{\gamma}\right\}$ & $\frac{3(r-1)(N-2)}{N}$ \\
\hline
$0$ & $1$ \\
\hline
\hline
\end{tabular}
\label{1:t22}
\end{table}

\section{Preliminary}


Denote by $C^{(N,r)}$ the subgroup of $\gf(r)^*$ generated by $\alpha^N$. Since $N|(m,q-1)$, the integer $(r-1)/(q-1)=q^{m-1}+q^{m-2}+\cdots+q+1$ is divisible by $N$, hence $\beta \in C^{(N,r)}$ and $\gf(q)^* \subset C^{(N,r)}$.

For any $u \in \gf(r)$, define
\begin{eqnarray} \label{2:eta} \eta_u^{(N,r)}=\sum_{z \in C^{(N,r)}}\psi(zu),\end{eqnarray}
where $\psi$ is the canonical additive character of $\gf(r)$, which is given by $\psi(x)=\exp\left(\frac{2 \pi i}{p} \tr_p(x)\right)$, here $\tr_p$ is the trace function from $\gf(r)$ to $\gf(p)$. Trivially $\eta_0^{(N,r)}=\frac{r-1}{N}$. If $u \ne 0$, the term $\eta_u^{(N,r)}$, which is called a ``Gaussian period'', depends only on the coset $uC^{(N,r)}$ in $\gf(r)^*$. There are $N$ such cosets, corresponding to $u=1,\alpha,\ldots,\alpha^{N-1}$ respectively, hence there are $N$ such Gaussian periods.

Recall from \cite[Lemma 5]{DL2} (see also \cite{DL1,WT}) that for any $(a,b) \in \gf(r)^2$, the Hamming weight of the codeword $\cc_{(a,b)}$ is given by
\begin{eqnarray} \label{2:weight} \omega\left(\cc_{(a,b)}\right)=\frac{h(r-1)}{q}-\lambda(a,b), \end{eqnarray}
where the ``modified'' weight $\lambda(a,b)$ is defined by
\[\lambda(a,b)=\frac{hN}{eq}\sum_{i=1}^e \eta_{(a+\beta^i b)g^i}^{(N,r)}. \]
It suffices to study $\lambda(a,b)$ only.

If $(a,b) \ne (0,0)$ while $a+\beta^t b=0$ for some $t$, $1 \le t \le e$, then $a=-\beta^tb$ and we have
\begin{eqnarray} \label{2:z2} \lambda\left(-\beta^t b,b\right)=\frac{hN}{eq} \left\{\frac{r-1}{N}+\sum_{\substack{i=1\\ i \ne t}}^e \eta_{bg^i(\beta^i-\beta^t)}^{(N,r)}\right\}. \end{eqnarray}
This is relatively easy to handle. On the other hand, if $\prod_{i=1}^e\left(a+\beta^i b\right) \ne 0$, this is more difficult. We have proved in \cite[Section 2]{X} the following: for any $c_1,\ldots,c_e \in \gf(r)^*$, let $\vec{c}=(c_1,\ldots,c_e)$ and define the set
\[\F(\vec{c})=\left\{(a,b) \in \gf(r)^2: \begin{array}{ll}
\left(a+\beta^i b\right)g^i c_i \in C^{(N,r)} \,\, \forall i \end{array}
\right\}.\]
Then for any $(a,b) \in \F(\vec{c})$, we have
\begin{eqnarray} \label{2:lambda}
\lambda(a,b)=\frac{hN}{eq}\sum_{i=1}^{e} \eta_{c_i^{-1}}^{(N,r)}.
\end{eqnarray}
Moreover,
\begin{eqnarray} \label{2:fc}
f(\vec{c}):=\#\F(\vec{c})=\frac{r-1}{N^e}\sum_{\substack{\chi_i^N=\epsilon\\
\chi_1, \ldots, \chi_{e-1}}} f_{\chi_1,\ldots,\chi_{e-1}}(\vec{c}),
\end{eqnarray}
where the sum is over all multiplicative characters $\chi_i$'s of $\gf(r)^*$ such that $\chi_i^N=\epsilon$, $\epsilon$ being the principal character, and
\[f_{\chi_1,\ldots,\chi_{e-1}}(\vec{c})=\sum_{b \in \gf(r)}\prod_{i=1}^{e-1} \chi_i\left(\xi_i(b+\mu_i)\right) ,\]
where \begin{eqnarray} \label{2:par}
\xi_i:=g^{i} (1-\beta^i) c_ic_e^{-1}, \quad \mu_i:=\frac{\beta^i}{1-\beta^i}, \quad i=1,2, \ldots, e-1.
\end{eqnarray}


\section{The case for $e=3$ and $N \ge 2$: the Jacobi sums}

It is apparent from (\ref{2:lambda}) and (\ref{2:fc}) that to determine the weight distribution of the code, we would need 1) the values of the Gaussian periods $\eta_u^{(N,r)}$, and 2) the values of the character sums $f(\vec{c})$ for each $\vec{c}$. Given the current status of knowledge, this is possible only for very special parameters $e$ and $N$. For $e=2$, the cases being treated in \cite{DL2,DL1}, the character sum $f(\vec{c})$ is directly related to counting the number of points on the curve $y^N=\xi_1(x+\mu_1),\,x,y \in \gf(r)$. This is a curve of genus zero. For cases such as $e=3, N=2$ (\cite{WT}), $e=4,N=2$ (\cite{X}) and $e=N=3$ (\cite{X2}), while the techniques are different, the main ideas can be summarized as reducing the calculation to counting the number of points on some elliptic curves. Here we observe that for $e=3$, in general, a more natural way to compute $f(\vec{c})$ is via the Jacobi sums.

Now we focus on the case that $e=3$ and $N \ge 2$. The parameters are
\[\beta=\alpha^{(r-1)/3}, \quad g=\alpha^{(q-1)/h}, \quad N=\gcd\left(m, \frac{3(q-1)}{h}\right),\quad 3|h \mbox{ and } h|(q-1). \]
Hence $\beta^3=1$ and $1+\beta+\beta^2=0$. Note that $\beta$ and any $a \in \gf(q)^*$ are both $N$-th powers in $\gf(r)$. Let us fix a multiplicative character $\chi$ of order $N$ defined on $\gf(r)^*$, then $\chi(\beta)=1$ and $\chi(a)=1$ for any $a \in \gf(r)^*$. In particular $\chi(-1)=1$. Let $\eta_{u}^{(N,r)}, u=1,\alpha,\ldots, \alpha^{N-1}$ be the $N$ Gaussian periods with respect to $C^{(N,r)}$.

Since all the characters of order dividing $N$ are given by $\chi,\chi^2,\ldots,\chi^N=\epsilon$, from (\ref{2:fc}), for $e=3$ we can write $f(\vec{c})$ as
\begin{eqnarray*}
f(\vec{c})&:=&\#\F(\vec{c})=\frac{r-1}{N^3} \sum_{b \in \gf(r)} \sum_{i=1}^N \chi^i\left(\xi_1(b+\mu_1)\right) \sum_{j=1}^N \chi^j\left(\xi_2(b+\mu_2)\right), \end{eqnarray*}
where $\xi_1,\xi_2,\mu_1,\mu_2$ are given in (\ref{2:par}), and
\begin{eqnarray} \label{3:lambda} \mu:=\mu_1-\mu_2=\frac{\beta}{1-\beta^2}. \end{eqnarray}
Changing variables $b \to \mu b$ and changing the order of summation, we find
\begin{eqnarray} \label{3:fc} f(\vec{c})=\frac{r-1}{N^3}\sum_{1 \le i,j \le N} \chi^i(\xi_1\mu) \chi^j(\xi_2\mu) \sum_{\substack{a+b=1\\
 a,b \in \gf(r)}} \chi^i\left(a\right) \chi^j\left(b\right).\end{eqnarray}
Here the inner sum is the so-called Jacobi sum defined by
\[J(\chi^i,\chi^j):=\sum_{\substack{a+b=1\\
 a,b \in \gf(r)}} \chi^i\left(a\right) \chi^j\left(b\right). \]
Using basic properties of Jacobi sums (\cite{IR}), we find that
\[J(\chi^N,\chi^N)= \sum_{\substack{a+b=1\\
0 \ne a,b \in \gf(r)}}1=r-2; \]
if $1 \le i,j \le N-1$ with $i+j=N$, then
\[J(\chi^i,\chi^j)=-\chi^i(-1)=-1;\]
if $i=N$, then
\[\sum_{j=1}^N\chi^j(\xi_2\mu) \sum_{\substack{a+b=1\\
a,b \in \gf(r)\\
a \ne 0}}\chi^j(b)=r-1-N\delta_N(\xi_2\mu); \]
if $j=N$, then
\[\sum_{i=1}^N\chi^i(\xi_1\mu) \sum_{\substack{a+b=1\\
a,b \in \gf(r)\\
b \ne 0}}\chi^i(a)=r-1-N\delta_N(\xi_1\mu). \]
Here we define $\delta_N(y):=1$ if $y \in C^{(N,r)}$ and $\delta_N(y):=0$ if $y \in \gf(r)-C^{(N,r)}$. We have also applied the orthogonal relations
\[\sum_{b \in \gf(r)}\chi^j(b)=\left\{\begin{array}{lll}
0&:& \mbox{ if } 1 \le j \le N-1,\\
r-1&:& \mbox{ if } j=N, \end{array}\right.\]
and
\[\sum_{j=1}^N\chi^j(a)=N\delta_N(a), \quad \forall \, a \in \gf(r)^*. \]
From these equations we can obtain
\begin{eqnarray} \label{3:nfc} f(\vec{c})=\frac{r-1}{N^3}\Bigl\{r+1-N\left(\delta_N(\xi_1\mu)+\delta_N(\xi_2 \mu)+\delta_N(\xi_1\xi_2^{-1})\right)+A\Bigr\},\end{eqnarray}
where
\[A=\sum_{\substack{1 \le i,j \le N-1\\
i+j \ne N}} \chi^i(\xi_1\mu) \chi^j(\xi_2\mu) J(\chi^i,\chi^j).\]
In order to evaluate the Jacobi sums $J(\chi^i,\chi^j)$, it is natural to consider the Gauss sums. For any multiplicative character $\eta$ of $\gf(r)$, the Gauss sum $\tau(\eta)$ is defined by
\[\tau(\eta)=\sum_{x \in \gf(r)} \eta(x) \psi(x). \]
Jacobi sums and Gauss sums are related (see \cite{IR}). In particular, for $1 \le i,j \le N-1$ with $i+j \ne N$, we have
\begin{eqnarray} \label{3:jacobi} J(\chi^i,\chi^j)=\frac{\tau(\chi^i)\tau(\chi^j)}{\tau(\chi^{i+j})}. \end{eqnarray}

\subsection{The case (1) of Theorem \ref{thm1}}

Under the assumptions of (1) in Theorem \ref{thm1}, from \cite[Proposition 20, Proposition 1]{MY} we derive easily that
\[\tau(\chi^i)=(-1)^i\sqrt{r}, \quad 1 \le i \le N-1,\]
Therefore $J(\chi^i,\chi^j)=\sqrt{r}$ for any $1 \le i,j \le N-1$ with $i+j \ne N$. We obtain
\[A=\sqrt{r}\sum_{\substack{1 \le i,j \le N-1\\
i+j \ne N}} \chi^i(\xi_1\mu) \chi^j(\xi_2\mu)=\sqrt{r}\left\{\left(\sum_{i=1}^{N-1}\chi^i(\xi_1\mu)\right)
\left(\sum_{j=1}^{N-1}\chi^j(\xi_2\mu)\right)-\sum_{i=1}^{N-1}\chi^i(\xi_1\xi_2^{-1})\right\}.\]
This again can be simplified as
\begin{eqnarray*} A&=&\sqrt{r}\left\{\left(N\delta_N(\xi_1\mu)-1\right)\left(N\delta_N(\xi_2\mu)-1\right)
-N\delta_N(\xi_1\xi_2^{-1})+1\right\}\\
&=&\sqrt{r} \left\{N^2\delta_N(\xi_1\mu) \delta_N(\xi_2\mu) -N(\delta_N(\xi_1\mu)+\delta_N(\xi_2\mu)+\delta_N(\xi_1\xi_2^{-1}))+2\right\}.  \end{eqnarray*}
Plugging this into (\ref{3:nfc}), and noting that $\beta \in C^{(N,r)}$, $1+\beta=-\beta^2 \in C^{(N,r)}$, we obtain the formula
\begin{eqnarray*}
f(\vec{c})&=&\frac{r-1}{N^3} \Bigl\{r+1-N\bigl\{\delta_N(gc_2c_1^{-1})+\delta_N(g^2c_2c_3^{-1})+\delta_N(gc_1c_3^{-1})\bigr\}+\sqrt{r} \times\\
&&  \bigl\{N^2\delta_N(g^2c_2c_3^{-1}) \delta_N(gc_1c_3^{-1}) -N(\delta_N(g^2c_2c_3^{-1})+\delta_N(gc_1c_3^{-1})+\delta_N(gc_2c_1^{-1}))+2\bigr\}\Bigr\}.
\end{eqnarray*}

\subsection{The case (2) of Theorem \ref{thm1}}

Under the assumptions of (2) in Theorem \ref{thm1}, from \cite[Proposition 20, Proposition1]{MY} we derive easily that
\[\tau(\chi^i)=(-1)^{\gamma+1}\sqrt{r}, \quad 1 \le i \le N-1,\]
where the integer $\gamma=\frac{sm}{2j}$ is defined in Theorem \ref{thm1}. So $J(\chi^i,\chi^j)=(-1)^{\gamma+1}\sqrt{r}$ for any $1 \le i,j \le N-1$ with $i+j \ne N$. Using this value we obtain a similar expression
\begin{eqnarray*}
f(\vec{c})&=&\frac{r-1}{N^3} \Bigl\{r+1-N\bigl\{\delta_N(gc_2c_1^{-1})+\delta_N(g^2c_2c_3^{-1})+\delta_N(gc_1c_3^{-1})\bigr\}-(-1)^{\gamma}\sqrt{r} \times\\
&&  \bigl\{N^2\delta_N(g^2c_2c_3^{-1}) \delta_N(gc_1c_3^{-1}) -N(\delta_N(g^2c_2c_3^{-1})+\delta_N(gc_1c_3^{-1})+\delta_N(gc_2c_1^{-1}))+2\bigr\}\Bigr\}.
\end{eqnarray*}

For other cases which are not covered by the assumptions of Theorem \ref{thm1}, using the Hasse-Davenport relation (\cite{IR}), it is possible to find all values of the Gauss sums $\tau(\chi^i)$, especially when $N$ is relatively small. However, even for $N=3$ (\cite{X2}), as we have seen, the results could be quite complicated, so are the formulas for $f(\vec{c})$.

\section{Proof of Theorem \ref{thm1}}

We will consider (2) of Theorem \ref{thm1} first. It is known from \cite[Proposition 20]{MY} that for this case we have
\begin{eqnarray*} \eta_1^{(N,r)}&=&\frac{-(-1)^{\gamma}(N-1)\sqrt{r}-1}{N}, \\ \eta_{\alpha^i}^{(N,r)}&=&\frac{(-1)^{\gamma}\sqrt{r}-1}{N}, \quad \forall \, 1 \le i \le N-1. \end{eqnarray*}
Let us define
\begin{eqnarray} \label{4:eta1} \eta_1:=\eta_1^{(N,r)}, \quad \eta_2:=\eta_{\alpha}^{(N,r)}.\end{eqnarray}
We summarize our argument as follows: from the equation (\ref{2:lambda}), the weight $\lambda(a,b)$ is a linear combination of $\eta_1,\eta_2$ with at most four distinct values, depending only on the cosets $c_iC^{(N,r)}$, $i=1,2,3$ such that $(a,b) \in \F(\vec{c})$. The size $f(\vec{c})=\#\F(\vec{c})$ has been obtained from the last section. For each $c_i$ there are $N$ cosets, which we represent simply as $c_i =1,\alpha, \alpha^2,\ldots,\alpha^{N-1}$. So we have all the ingredients to compute the weight distribution. The final results will depend on whether or not $N|\frac{q-1}{h}$.

\subsection{Proof of (2.1) in Theorem \ref{thm1}}
If $N|\frac{q-1}{h}$, then $g \in C^{(N,r)}$ is an $N$-th power. For each $\vec{c}=(c_1,c_2,c_3) \in \{1,\alpha,\ldots,\alpha^{N-1}\}^{3}$, we find from the end of the last section
\begin{eqnarray*}
f(\vec{c})&=&\frac{r-1}{N^3} \Bigl\{r+1-N\bigl\{\delta_N(c_2c_1^{-1})+\delta_N(c_2c_3^{-1})+\delta_N(c_1c_3^{-1})\bigr\}-(-1)^{\gamma}\sqrt{r} \times\\
&&  \bigl\{N^2\delta_N(c_2c_3^{-1}) \delta_N(c_1c_3^{-1}) -N(\delta_N(c_2c_3^{-1})+\delta_N(c_1c_3^{-1})+\delta_N(c_2c_1^{-1}))+2\bigr\}\Bigr\},
\end{eqnarray*}
and for each $(a,b)\in \F(\vec{c})$,
\[\lambda(a,b)=\frac{hN}{3q}\sum_{i=1}^3\eta_{c_i^{-1}}^{(N,r)}. \]

If $c_1=c_2=c_3=1$, i.e., $c_i \in C^{(N,r)} \, \forall i$, the weight is $\lambda(a,b)=\frac{hN\eta_1}{q}$, and the number of such $(a,b)$'s counted in this $\F(\vec{c})$ is
\[f(\vec{c})=\frac{r-1}{N^3}\Bigl\{r+1-3N-(-1)^{\gamma}\sqrt{r}\left(N^2-3N+2\right)\Bigr\}.\]

If one of $c_1,c_2,c_3$ is $1$ and the other two are not $1$ (for example, $c_1=1,c_2=\alpha,c_3=\alpha^2$), the weight is $\lambda(a,b)=\frac{hN}{3q}(\eta_1+2\eta_2)$. Denote by $\A_1$ the set of all the $\vec{c}=(c_1,c_2,c_3)$'s with this property. Clearly $\#\A_1=3(N-1)^2$. For $\vec{c} \in \A_1$, if $c_1=1$, then $c_2 \ne 1$, hence $\delta_N(c_2c_1^{-1})=0$; if $c_2=1$, then $c_1 \ne 1$, so again $\delta_N(c_2c_1^{-1})=0$; if $c_3=1$, then $c_1,c_2 \ne 1$, then $c_2c_1^{-1} \in C^{(N,r)}$ so that $\delta_N(c_2c_1^{-1})=1$ if and only if $c_1=c_2=\alpha^i, 1 \le i \le N-1$. Thus we have proved that
\[\sum_{\vec{c} \in \A_1}\delta_N(c_2c_1^{-1})=N-1. \]
Similarly
\[\sum_{\vec{c} \in \A_1}\delta_N(c_2c_3^{-1})=\sum_{\vec{c} \in \A_1}\delta_N(c_1c_3^{-1})=N-1. \]
It is also easy to check that
\[\sum_{\vec{c} \in \A_1}\delta_N(c_2c_3^{-1})\delta_N(c_1c_3^{-1})=0. \]
So we obtain
\begin{eqnarray*}\sum_{\vec{c} \in \A_1}f(\vec{c})&=&\frac{r-1}{N^3}\left\{3(N-1)^2(r+1)-3N(N-1)-(-1)^{\gamma}\sqrt{r}\left(-3N(N-1)+6(N-1)^2\right)\right\}\\
&=&\frac{3(r-1)(N-1)}{N^3}\bigl\{r(N-1)-\sqrt{r}(-1)^{\gamma}(N-2)-1\bigr\}.
\end{eqnarray*}

If two of $c_1,c_2,c_3$ are equal to $1$ but the third is not equal to $1$ (for example, $c_1=c_2=1,c_3=\alpha$), the weight is $\lambda(a,b)=\frac{hN}{3q}(2\eta_1+\eta_2)$. Denote by $\A_2$ the set of all the $\vec{c}=(c_1,c_2,c_3)$'s with this property. Clearly $\#\A_2=3(N-1)$. Using similar arguments as above, we find
\[\sum_{\vec{c} \in \A_2}\delta_N(c_2c_1^{-1})=\sum_{\vec{c} \in \A_2}\delta_N(c_2c_3^{-1})=\sum_{\vec{c} \in \A_1}\delta_N(c_1c_3^{-1})=N-1, \]
and
\[\sum_{\vec{c} \in \A_2}\delta_N(c_2c_3^{-1})\delta_N(c_1c_3^{-1})=0. \]
So we obtain
\begin{eqnarray*}\sum_{\vec{c} \in \A_2}f(\vec{c})&=&\frac{r-1}{N^3}\left\{3(N-1)(r+1)-3N(N-1)-(-1)^{\gamma}\sqrt{r}\left(-3N(N-1)+6(N-1)\right)\right\}\\
&=&\frac{3(r-1)(N-1)}{N^3}\bigl\{r+\sqrt{r}(-1)^{\gamma}(N-2)-N+1\bigr\}.
\end{eqnarray*}

If none of $c_1,c_2,c_3$ is equal to $1$, the weight is $\lambda(a,b)=\frac{hN\eta_2}{q}$. Denote by $\A_3$ the set of all the $\vec{c}=(c_1,c_2,c_3)$'s with this property. Clearly $\#\A_3=(N-1)^3$. It is easy to check that
\[\sum_{\vec{c} \in \A_2}\delta_N(c_2c_1^{-1})=\sum_{\vec{c} \in \A_2}\delta_N(c_2c_3^{-1})=\sum_{\vec{c} \in \A_1}\delta_N(c_1c_3^{-1})=(N-1)^2, \]
and
\[\sum_{\vec{c} \in \A_2}\delta_N(c_2c_3^{-1})\delta_N(c_1c_3^{-1})=N-1. \]
So we obtain
\begin{eqnarray*}\sum_{\vec{c} \in \A_2}f(\vec{c})&=&\frac{(r-1)(N-1)}{N^3}\bigl\{r(N-1)^2+\sqrt{r}(-1)^{\gamma}(N-2)-(N-1)(2N+1)\bigr\}.
\end{eqnarray*}

Finally, we need to deal with the easy case that $(a,b) \ne (0,0)$ but $a+\beta^tb=0$ for some $t, 1 \le t \le 3$. First, it can be seen that $\beta^i-\beta^j \in (1-\beta)C^{(N,r)}$ for any $1 \le i \ne j \le 3$. Since $g \in C^{(N,r)}$, if $t=1$, then $a=-\beta b$, from (\ref{2:z2}) we find that
\[\lambda(-\beta b,b)=\frac{hN}{3q}\left\{\frac{r-1}{N}+2 \eta_{b(1-\beta)}^{(N,r)}\right\}. \]
As $b$ varies in $\gf(r)^*$, more precisely, for any $b \in (1-\beta)^{-1}C^{(N,r)}$, we have
\[\lambda(-\beta b,b)=\frac{hN}{3q}\left\{\frac{r-1}{N}+2 \eta_1\right\}, \]
the number of such $b$'s is $(r-1)/N$. For any $b \in \gf(r)^*-(1-\beta)^{-1}C^{(N,r)}$, the weight is
\[\lambda(-\beta b,b)=\frac{hN}{3q}\left\{\frac{r-1}{N}+2 \eta_2\right\}, \]
the number of such $b$'s is $(r-1)(N-1)/N$. For $t=2,3$, the results are the same. We summarize the results in Table \ref{4:t1}.

\begin{table}[ht] 
\caption{The case for (2.1) of Theorem \ref{thm1}}
\centering
\begin{tabular}{|c|| c|}
\hline
\hline
Weight $\lambda(a,b)$& Frequency  \\
[0.5ex]
\hline
$\frac{hN}{q}\eta_1$ & $\frac{r-1}{N^3}\bigl\{r-\sqrt{r}(-1)^{\gamma}(N^2-3N+2)-3N+1\bigr\}$ \\
\hline
$\frac{hN}{q}\eta_2$ & $\frac{(r-1)(N-1)}{N^3}\bigl\{r(N-1)^2+\sqrt{r}(-1)^{\gamma}(N-2)-(N-1)(2N+1)\bigr\}$ \\
\hline
$\frac{hN}{3q}\left\{\eta_1+2\eta_2\right\}$ & $\frac{3(r-1)(N-1)}{N^3}\bigl\{r(N-1)-\sqrt{r}(-1)^{\gamma}(N-2)-1\bigr\}$ \\
\hline
$\frac{hN}{3q}\left\{2\eta_1+\eta_2\right\}$ &  $\frac{3(r-1)(N-1)}{N^3}\bigl\{r+\sqrt{r}(-1)^{\gamma}(N-2)-N+1\bigr\}$ \\
\hline
$\frac{hN}{3q}\left\{\frac{r-1}{N}+2\eta_1\right\}$ & $\frac{3(r-1)}{N}$ \\
\hline
$\frac{hN}{3q}\left\{\frac{r-1}{N}+2\eta_1\right\}$ & $\frac{3(r-1)(N-1)}{N}$ \\
\hline
\hline
\end{tabular}
\label{4:t1}
\end{table}
Using the values $\eta_1,\eta_2$ from (\ref{4:eta1}) and noting the relation between the Hamming weight $w(\cc(a,b))$ and $\lambda(a,b)$ in (\ref{2:weight}), we obtain Table \ref{1:t21}. The extra Hamming weight $0$ with frequency $1$ in Table \ref{1:t21} comes from the codeword with $a=b=0$. This finishes the proof of (2.1) in Theorem \ref{thm1}. $\square$

\subsection{Proof of (2.2) in Theorem \ref{thm1}}

If $N \dagger \frac{q-1}{h}$, then $g \notin C^{(N,r)}, g^2 \notin C^{(N,r)},g^3 \in C^{(N,r)}$. For each $\vec{c}=(c_1,c_2,c_3) \in \{1,\alpha,\ldots,\alpha^{N-1}\}^{3}$, we find by the end of the last section
\begin{eqnarray*}
f(\vec{c})& =&\frac{r-1}{N^3} \Bigl\{r+1-N\bigl\{\delta_N(gc_2c_1^{-1})+\delta_N(g^2c_2c_3^{-1})+\delta_N(gc_1c_3^{-1})\bigr\}-(-1)^{\gamma}\sqrt{r} \times\\
&&\bigl\{N^2\delta_N(g^2c_2c_3^{-1}) \delta_N(gc_1c_3^{-1}) -N\bigl\{\delta_N(g^2c_2c_3^{-1})+\delta_N(gc_1c_3^{-1})+\delta_N(gc_2c_1^{-1})\bigr\}+2\bigr\}\Bigr\},
\end{eqnarray*}
and for each $(a,b)\in \F(\vec{c})$,
\[\lambda(a,b)=\frac{hN}{3q}\sum_{i=1}^3\eta_{c_i^{-1}}^{(N,r)}. \]

The argument is very similar. However, for the sake of completeness, we give a full account of the proof here. If $c_1=c_2=c_3=1$, i.e., $c_i \in C^{(N,r)} \, \forall i$, then the weight is $\lambda(a,b)=\frac{hN\eta_1}{q}$, and the number of such $(a,b)$'s counted in this $\F(\vec{c})$ is
\[f(\vec{c})=\frac{r-1}{N^3}\Bigl\{r-2\sqrt{r}(-1)^{\gamma}+1\Bigr\}.\]
If one of $c_1,c_2,c_3$ is $1$ and the other two are not $1$, then the weight is $\lambda(a,b)=\frac{hN}{3q}(\eta_1+2\eta_2)$. Denote by $\A_1$ the set of all the $\vec{c}=(c_1,c_2,c_3)$'s with this property, whose size is clearly $3(N-1)^2$. It is easy to check that
\[\sum_{\vec{c} \in \A_1}\delta_N(gc_2c_1^{-1})=\sum_{\vec{c} \in \A_1}\delta_N(g^2c_2c_3^{-1})=\sum_{\vec{c} \in \A_1}\delta_N(gc_1c_3^{-1})=3N-4, \]
and
\[\sum_{\vec{c} \in \A_1}\delta_N(g^2c_2c_3^{-1})\delta_N(gc_1c_3^{-1})=3. \]
So we obtain
\begin{eqnarray*}\sum_{\vec{c} \in \A_1}f(\vec{c})&=&\frac{3(r-1)}{N^3}\bigl\{r(N-1)^2-2\sqrt{r}(-1)^{\gamma}-2N^2+2N+1\bigr\}.
\end{eqnarray*}
If two of $c_1,c_2,c_3$ are equal to $1$ but the third is not equal to $1$, then the weight is $\lambda(a,b)=\frac{hN}{3q}(2\eta_1+\eta_2)$. Denote by $\A_2$ the set of all the $\vec{c}=(c_1,c_2,c_3)$'s with this property, whose size is clearly $3(N-1)$. It is easy to check that
\[\sum_{\vec{c} \in \A_2}\delta_N(gc_2c_1^{-1})=\sum_{\vec{c} \in \A_2}\delta_N(g^2c_2c_3^{-1})=\sum_{\vec{c} \in \A_1}\delta_N(gc_1c_3^{-1})=2, \]
and
\[\sum_{\vec{c} \in \A_2}\delta_N(g^2c_2c_3^{-1})\delta_N(gc_1c_3^{-1})=0. \]
So we obtain
\begin{eqnarray*}\sum_{\vec{c} \in \A_2}f(\vec{c})&=&\frac{3(r-1)}{N^3}\bigl\{r(N-1)+2\sqrt{r}(-1)^{\gamma}-N-1\bigr\}.
\end{eqnarray*}
If none of $c_1,c_2,c_3$ is equal to $1$, then the weight is $\lambda(a,b)=\frac{hN\eta_2}{q}$. Denote by $\A_3$ the set of all the $\vec{c}=(c_1,c_2,c_3)$'s with this property, whose side is clearly $(N-1)^3$. It is easy to check that
\[\sum_{\vec{c} \in \A_2}\delta_N(gc_2c_1^{-1})=\sum_{\vec{c} \in \A_2}\delta_N(g^2c_2c_3^{-1})=\sum_{\vec{c} \in \A_1}\delta_N(gc_1c_3^{-1})=(N-1)(N-2), \]
and
\[\sum_{\vec{c} \in \A_2}\delta_N(g^2c_2c_3^{-1})\delta_N(gc_1c_3^{-1})=N-3. \]
So we obtain
\begin{eqnarray*}\sum_{\vec{c} \in \A_2}f(\vec{c})&=&\frac{r-1}{N^3}\bigl\{r(N-1)^3+2\sqrt{r}(-1)^{\gamma}-(N-1)(2N^2-4N-1)\bigr\}.
\end{eqnarray*}
We also need to deal with the easy case that $(a,b) \ne (0,0)$ but $a+\beta^tb=0$ for some $t, 1 \le t \le 3$. Since $\beta^i-\beta^j \in (1-\beta)C^{(N,r)}$ for any $1 \le i \ne j \le 3$. Since $g,g^2 \notin C^{(N,r)}, g^3 \in C^{(N,r)}$, if $t=1$, then $a=-\beta b$, from (\ref{2:z2}) we find that
\[\lambda(-\beta b,b)=\frac{hN}{3q}\left\{\frac{r-1}{N}+ \eta_{bg^2(1-\beta)}^{(N,r)}+\eta_{b(1-\beta)}^{(N,r)}\right\}. \]
As $b$ varies in $\gf(r)^*$, we find that the weight is
\[\lambda(-\beta b,b)=\frac{hN}{3q}\left\{\frac{r-1}{N}+ \eta_1+\eta_2\right\} \]
for $2(r-1)/N$ many such $b$'s, and the weight is
\[\lambda(-\beta b,b)=\frac{hN}{3q}\left\{\frac{r-1}{N}+2 \eta_2\right\} \]
for $(r-1)(N-2)/N$ many such $b$'s. For $t=2,3$, the results are the same.

Summarizing the above results, using the values $\eta_1,\eta_2$ from (\ref{4:eta1}) and the relation between the Hamming weight $w(\cc(a,b))$ and $\lambda(a,b)$ in (\ref{2:weight}), we obtain Table \ref{1:t21}. This finishes the proof of (2.2) in Theorem \ref{thm1}. $\square$

\subsection{Proofs of (1.1) and (1.2) in Theorem \ref{thm1}}

Now we consider (1) of Theorem \ref{thm1}. It is known from \cite[Proposition 20]{MY} that for this case $N$ is even and
\begin{eqnarray*} \eta_{\alpha^{N/2}}^{(N,r)}&=&\frac{(N-1)\sqrt{r}-1}{N}, \\ \eta_{\alpha^i}^{(N,r)}&=&\frac{-\sqrt{r}-1}{N}, \quad \forall \, 0 \le i \le N-1, i \ne N/2. \end{eqnarray*}
Let us define
\begin{eqnarray} \label{4:eta2} \eta_1:=\eta_{\alpha^{N/2}}^{(N,r)}, \quad \eta_2:=\eta_{1}^{(N,r)}.\end{eqnarray}
The arguments are very similar, given the expression of $f(\vec{c})$ in the last section, so we omit the details. In summary, after considering the cases 1) $c_1=c_2=c_3=\alpha^{N/2}$, 2) exactly one of $c_1,c_2,c_3$ is $\alpha^{N/2}$, 3) exactly two of $c_1,c_2,c_3$ are $\alpha^{N/2}$, 4) none of $c_1,c_2,c_3$ is $\alpha^{N/2}$, and 5) the cases such that $(a,b) \ne (0,0)$ but $a+\beta^tb=0$ for some $1 \le t \le 3$, we can obtain Table \ref{1:t11} if $N|\frac{q-1}{h}$ and Table \ref{1:t12} if $N \dagger \frac{q-1}{h}$. This finishes the proof of (1) in Theorem \ref{thm1}. Now the proof of Theorem \ref{thm1} is complete. $\square$

\end{document}